\documentclass{amsart}

\usepackage{amsthm,amsmath,amssymb,amscd,amsfonts}

\theoremstyle{plain}
\newtheorem{theorem}{Theorem}[section]

\newtheorem{proposition}[theorem]{Proposition}
\newtheorem{corollary}[theorem]{Corollary}

\theoremstyle{definition}

\newtheorem{conjecture}[theorem]{Conjecture}
\newtheorem*{acknowledgement}{Acknowledgement}

\theoremstyle{remark}
\newtheorem{remark}[theorem]{Remark}

\renewcommand{\P}{{\mathbb{P}}}
\DeclareMathOperator{\Aut}{{Aut}}

\DeclareMathOperator{\Chow}{{Chow}}

\DeclareMathOperator{\oka}{{\mathcal O}}

\begin{document}

\title[Effective finiteness theorems for maps]{Effective finiteness theorems for maps between canonically polarized compact complex manifolds}

\author[G. Heier]{Gordon Heier}

\address{Ruhr-Universit\"at Bochum\\
Fakult\"at f\"ur Mathematik\\
D-44780 Bochum\\
Germany}

\email{heier@cplx.ruhr-uni-bochum.de}

\subjclass[2000]{14Q20, 32H02, 32H35, 32J18}

\begin{abstract}
Effective bounds for the finite number of surjective holomorphic maps between canonically polarized compact complex manifolds of any dimension with fixed domain are proven. Both the case of a fixed target and the case of varying targets are treated. In the case of varying targets, bounds on the complexity of Chow varieties are used.
\end{abstract}

\maketitle

\section{Effective bounds for automorphism groups}\label{autsection}
Hurwitz proved the following effective finiteness theorem on Riemann surfaces.
\begin{theorem}[\cite{Hur}]\label{Hurbound}
Let $X$ be a smooth compact complex curve of genus $g\geq 2$. Then
the group $\Aut(X)$ of holomorphic automorphisms of $X$ satisfies
\begin{equation*}
  \#\Aut(X)\leq84(g-1).
\end{equation*}
\end{theorem}
For many years after Hurwitz's proof, this bound has been known to be sharp only for $g=3$ and $g=7$, in which cases there exist, respectively, the classical examples of the Klein quartic in $\P^2$ given by the homogeneous equation $X^3Y+Y^3Z+Z^3X=0$ and the Fricke curve with automorphism group ${\rm PSL}(2,8)$. Using the theory of finite groups, it was established only in the 1960's by Macbeath that there are infinitely many $g$ for which the above bound is sharp (see \cite{Macbeath}).\par
Xiao was able to establish the following direct (and clearly sharp due to the above) generalization of Hurwitz's theorem.
\begin{theorem}[\cite{Xi1}, \cite{Xi2}]
Let $X$ be a $2$-dimensional minimal compact complex
manifold of general type. Then
\begin{equation*}
  \#\Aut(X)\leq(42)^2K_X^2.
\end{equation*}
\end{theorem}
In arbitrary dimension, the automorphism group of a smooth compact complex manifold of general type is still known to be finite because of the finiteness theorem of Kobayashi-Ochiai (\cite{KobOch}), which we shall state in the next section. One is tempted to conjecture that in the case of the canonical line bundle being big and nef or even ample, there is an upper bound of the form $C_nK_X^n$. The preprint \cite{Ts} makes an attempt to prove this conjecture.\par
In the paper \cite{Sza}, Szab\'o was able to establish the following effective polynomial upper bound in arbitrary dimension.
\begin{theorem}[\cite{Sza}]\label{szabobound}
Let $X$ be an $n$-dimensional compact complex manifold whose canonical line bundle is big and nef. Then the number of birational automorphisms of $X$ is no more than
\begin{equation*}
(2(n+1)(n+2)!(n+2)K_X^n)^{16n3^n}.
\end{equation*}
\end{theorem}
The multiple $2(n+1)(n+2)!(n+2)K_X$ is large enough to give a birational morphism from $X$ to projective space. This is proven in \cite[page 8]{CaSchn}, using results of Demailly \cite{DBound} and Koll\'ar \cite{Koeffbase} on effective base point freeness of adjoint line bundles. The goal of \cite{CaSchn} is to obtain a polynomial bound for the special case of automorphism groups that are abelian.\par
In arbitrary dimension, effective pluricanonical (birational) embeddings are essential in proving finiteness statements of the type considered in this paper. They enable us to bring the problem into the context of projective varieties and to establish uniform boundedness. In the case of $K_X$ being ample, the following effective theorem on pluricanonical embeddings is available.
\begin{theorem}[\cite{heierefffreeness}, Corollary 1.14]\label{effpluri} If $X$ is a compact complex manifold of complex
dimension $n$ whose canonical line bundle $K_X$ is ample, then $mK_X$
is very ample for any integer
\begin{equation}\label{pluricaneffbound}
m\geq (e+\frac 1 2)n^\frac 7 3+\frac 1 2 n^\frac 5 3 + (e+\frac 1
2)n^\frac 4 3 + 3n+ \frac 1 2 n^\frac 2 3+5,
\end{equation}
where $e \approx 2.718$ is Euler's number.
\end{theorem}
From now on, we will set $k=k(n)$ to be the round-up of the effective very ampleness bound given in \eqref{pluricaneffbound}.\par
To our knowledge, Szab\'o's theorem is the one that provides the best bound at this point. However, its proof relies on several methods previously introduced by other authors (e.g., see \cite{HuckSauer}) and uses the classification of finite simple groups in an essential way. In light of this, the much more straightforward method of Howard and Sommese, which was introduced in \cite{HoSo}, still deserves to be noticed. Their method is actually not only applicable to automorphisms (see next section), and it represents an instance of a proof based entirely on boundedness and rigidity, which, technically speaking, is the main focus of the present paper. \par
Howard and Sommese prove for the case of a canonically polarized manifold that $\#Aut(X)$ is bounded from above by a number which depends only on the Chern numbers of $X$. Based on their result, we now state the following effective
finiteness theorem.
\begin{theorem}\label{HoSoAutbound}
Let $X$ be a compact complex manifold of dimension $n$ whose
canonical line bundle is ample. Then
\begin{equation*}
\#\Aut(X) \leq \left((n+1)^2k^nn!2^{n^2}(2k)^{\frac 1 2 n(n+1)}(1+2kn)^nK_X^n\right)^{((k^nK_X^n+n)^2-1} .
\end{equation*}
\end{theorem}
Before we prove this theorem, we need to prove two auxiliary propositions which make the method of Howard and Sommese entirely effective. The first proposition will be used to bound the dimension of the target projective space for the
pluricanonical embedding given by $kK_X$. It is a standard argument.
\begin{proposition}\label{proph0bound}
Let $X$ be an $n$-dimensional compact complex manifold and $L$ a very ample line bundle on $X$. Then
\end{proposition}
\begin{equation*}
h^0(X,L)\leq L^n+n.
\end{equation*}
\begin{proof}
We proceed by induction.\par
The case $n=1$ follows immediately from the Riemann-Roch Theorem.\par
Let $D$ be an effective divisor on $X$ such that $\oka_X(D)=L$.
One has the standard short exact sequence
\begin{equation*}
0\to\oka_X\to\oka_X(L)\to\oka_D(L)\to 0.
\end{equation*}
From this exact sequence, we obtain
\begin{equation*}
h^0(X,L)\leq h^0(X,\oka_X)+h^0(D,\oka_D(L)).
\end{equation*}
By induction, we find that
\begin{equation*}
h^0(X,L)\leq 1+(L_{|D})^{n-1}+n-1 = L^n+n.
\end{equation*}
\end{proof}
Secondly, we use a result of Demailly, Peternell and Schneider in \cite{DPS} to compute a bound for the Chern class intersection numbers that
occur in the well-known formula for the degree of the
($1$-codimensional) dual of a projective variety. Our effective result is the following.
\begin{proposition}\label{chernintersection}
Let $X$ be a compact complex manifold of dimension $n$ whose
canonical line bundle is ample. Let $k$ denote again the round-up of the constant defined in
\eqref{pluricaneffbound} Then the following holds for $i=1,\ldots,n$.
\begin{equation*}
|c_i(\Omega_X^1).K_X^{n-i}|\leq i!2^{in}(2k)^{\frac 1 2 i(i+1)}(1+2kn)^iK_X^n.
\end{equation*}
\end{proposition}
\begin{proof}
Recall that $k$ is such that $kK_X$ is very ample. It follows from the Castelnuovo-Mumford theory of
regularity that $\Omega_X^1(2kK_X)$ is generated by global sections and
therefore nef. We may thus apply \cite[Corollary 2.6]{DPS} to obtain
\begin{eqnarray}
  0&\leq&c_i(\Omega_X^1(2kK_X))K_X^{n-i}\nonumber\\
 &\leq& (c_1(\Omega_X^1(2kK_X)))^iK_X^{n-i}\nonumber\\
 &=&(c_1(\Omega_X^1)+2knK_X)^iK_X^{n-i}\nonumber\\
  &=&(1+2kn)^iK_X^n\label{B}
\end{eqnarray}
for $i=1,\ldots,n$.\par
In \cite[page 56]{Ful}, one finds the formula
\begin{equation}\label{chernclassformula}
c_i(\Omega_X^1(2kK_X))=\sum_{\nu=0}^{i}{n-\nu\choose i-\nu}c_\nu(\Omega_X^1)(2kK_X)^{i-\nu},
\end{equation}
which enables us to prove the Proposition by an induction.\par The inequality clearly holds in the case $i=1$.\par
For $1<i\leq n$, note that it follows from \eqref{chernclassformula} that 
\begin{eqnarray*}
c_i(\Omega_X^1(2kK_X))K_X^{n-i}&=&\left(\sum_{\nu=0}^{i}{n-\nu\choose i-\nu}c_\nu(\Omega_X^1)(2kK_X)^{i-\nu}\right)K_X^{n-i}\\
&=&\sum_{\nu=0}^{i}{n-\nu\choose i-\nu}c_\nu(\Omega_X^1)(2k)^{i-\nu}K_X^{n-\nu}.
\end{eqnarray*}
Taking absolute values, the triangle inequality yields
\begin{eqnarray*}
&&|c_i(\Omega_X^1)K_X^{n-i}|\\
&\leq&c_i(\Omega_X^1(2kK_X))K_X^{n-i}+\sum_{\nu=0}^{i-1}{n-\nu\choose i-\nu}|c_\nu(\Omega_X^1)(2k)^{i-\nu}K_X^{n-\nu}|\\
&\stackrel{\eqref{B}}{\leq}&(1+2kn)^iK_X^n+\sum_{\nu=0}^{i-1}{n-\nu\choose i-\nu}(2k)^{i-\nu}|c_\nu(\Omega_X^1)K_X^{n-\nu}|\\
&\stackrel{Ind.}{\leq}&(1+2kn)^iK_X^n+\sum_{\nu=0}^{i-1}{n-\nu\choose i-\nu}(2k)^{i-\nu}\nu! 2^{\nu n} (2k)^{\frac 1 2 \nu(\nu+1)}(1+2kn)^\nu K_X^n\\
&\leq&(1+2kn)^iK_X^n+i2^n(2k)^i(i-1)!2^{(i-1)n}(2k)^{\frac 1 2 i(i-1)}(1+2kn)^{i-1}K_X^{n}\\
&\leq&(1+2kn)^iK_X^n+2^n(2k)^{\frac 1 2 i(i+1)}i!2^{(i-1)n}(1+2kn)^{i-1}K_X^{n}\\
&\leq&(1+2kn)^iK_X^n+(2k)^{\frac 1 2 i(i+1)}i!2^{in}(1+2kn)^{i-1}K_X^{n}\\
&\leq&i!2^{in}(2k)^{\frac 1 2 i(i+1)}(1+2kn)^{i}K_X^{n}\quad \text{q.e.d.}
\end{eqnarray*}
\end{proof}
Now that we have all necessary effective tools at our disposal, we can proceed to the
\begin{proof}[Proof of Theorem \ref{HoSoAutbound}]
The proof given in \cite{HoSo} yields that
\begin{equation}\label{HoSoChern}
 \#\Aut(X)\leq \left(\sum_{\nu=0}^{n}(-1)^j(n+1-j)(kK_X)^{n-j}c_j(X)\right)^{(h^0(kK_X))^2-1}.
\end{equation}
Substituting the numerical bounds derived in Propositions \ref{proph0bound} and \ref{chernintersection} and estimating in an obvious way, we obtain that
\begin{equation*}
\#\Aut(X) \leq\left((n+1)^2k^nn!2^{n^2}(2k)^{\frac 1 2 n(n+1)}(1+2kn)^nK_X^n\right)^{(k^nK_X^n+n)^2-1}.
\end{equation*}
\end{proof}
\section{Effective finiteness theorems for maps with a fixed target}\label{fixed target}
For surjective meromorphic maps between compact complex spaces there is the following finiteness theorem due to Kobayashi-Ochiai.
\begin{theorem}[\cite{KobOch}]
Let $X$ be any compact complex space and $Y$ a compact complex space of general type. Then the number of surjective meromorphic maps between $X$ and $Y$ is finite.
\end{theorem}
There are no known effective versions of this theorem due to the fact that there
are no effective birational embedding theorems for manifolds with merely big
canonical line bundle in higher dimensions. The case of $X$ and $Y$ being smooth compact complex curves has been known for a long time as the Theorem of de Franchis, based on \cite{deF}. Not surprisingly, there are effective bounds in this case that depend only on the genus $g$ of $X$. However, these bounds are often obtained in the more general case of varying targets (see \cite{HS}, \cite{Guerra}) or in complete analogy to the higher dimensional case. Somewhat surprisingly, those authors that consider specifically the case of two {\it fixed} Riemann surfaces and investigate e.g. the induced homomorphisms on the first homology groups (as e.g. in \cite{tanabe}) do not seem to be able to do much better numerically than those who consider more general situations. All bounds are exponential in $g$, and the question of the true nature of the dependence on $g$ seems to be completely open. Since maps between fixed Riemann surfaces seem to be closer in spirit to automorphisms than to the case of varying targets, where the bound is not polynomial (see next section), and based on some other preliminary evidence, we venture the following conjecture.
\begin{conjecture}\label{deFranchisConj}
There is a polynomial function $B(g)$ with the following property. For two fixed smooth compact complex curves $X$ and $Y$ of genus at least $2$ with the genus of $X$ equal to $g$, and the number of surjective holomorphic maps from $X$ to $Y$ is no more than $B(g)$.
\end{conjecture}
As was already indicated in the previous section, the method of Howard and Sommese for automorphism groups can also be used to obtain a bound for the number of maps between any two fixed canonically polarized manifolds. The details of this straightforward generalization can be found in \cite{BD}. In fact, the bound one arrives at is the same as the expression we already encountered in \eqref{HoSoChern}. So we simply state the following theorem.
\begin{theorem}\label{efffixedtarget}
Let $X$ and $Y$ be fixed compact complex manifolds with ample canonical line bundles. Let $n$ be the dimension of $X$. Then the number of surjective holomorphic maps between $X$ and $Y$ is no more than 
\begin{equation*}
\left((n+1)^2k^nn!2^{n^2}(2k)^{\frac 1 2 n(n+1)}(1+2kn)^nK_X^n\right)^{(k^nK_X^n+n)^2-1} .
\end{equation*}
\end{theorem}
As we move on, we remark that a bound for the above theorem can also be obtained by using the Chow variety method discussed in the next section, since the graph of a surjective holomorphic map $X\to Y$ corresponds to an isolated point in a certain Chow variety of $X\times Y$. However, since this leads in fact to a worse bound, we will not discuss this in detail.
\section{Effective finiteness theorems for maps with varying targets}
The following theorem is often referred to as the Theorem of de Franchis-Severi. Its statement is obtained from the statement of the de Franchis Theorem by allowing the targets $Y$ to vary among smooth compact complex curves of genus at least $2$.
\begin{theorem}\label{deFSev}
Let $X$ be a smooth compact complex curve. Then the set of all holomorphic
maps $f:X\to Y$, where $Y$ is any (variable) smooth compact complex curve of genus
at least $2$, is finite.
\end{theorem}
In \cite{HS}, Howard and Sommese proved that if $X$ is of genus
$g$, the number of holomorphic maps in Theorem \ref{deFSev} modulo
automorphisms of the target spaces is no more than
\begin{equation}\label{HSbound}
\left(\frac 1 2
(2\sqrt{6}(g-1)+1)^{2+2g^2}g^2(g-1)(\sqrt{2})^{g(g-1)}+1\right).
\end{equation}
We denote this expression by $\mathcal S'(g)$. Since the
cardinality of the automorphism group of any one of the targets is at most
$\tfrac{1}{2}\cdot 84(g-1)$ due to Hurwitz, one can alternatively say that the
number of holomorphic maps in Theorem \ref{deFSev} is no more than
\begin{equation*}
  \mathcal S(g):= 42(g-1)\cdot\mathcal S'(g).
\end{equation*}
\begin{remark}
In their paper, Howard and Sommese apparently overlooked the fact
that their technique counts maps only modulo automorphisms of the targets. This
fact was observed by Kani in his paper \cite{Kani}. This oversight can, of course, easily be remedied by adding the factor $42(g-1)$. In Kani's paper, isomorphism classes of targets are counted instead of maps by means of a ``packing argument''. We chose to quote the result of Howard and Sommese because its proof is closer to the point of view taken in the present paper.
\end{remark}
It is certainly interesting to note that \cite[\S 4]{Kani} exhibits a relatively straightforward example of a series of Riemann surfaces $X$ that shows that the cardinality of the sets of maps defined in Theorem \ref{deFSev} cannot be bounded by a polynomial in $g$. In fact, what is shown is that the number of isomorphism classes of targets in these sets cannot be bounded by a polynomial. Therefore, Kani's example does not contradict our Conjecture \ref{deFranchisConj}. We would like to remark that the statements in \cite[page 802]{BD} and \cite[page 110]{TsaiIMRN}, which say that there is such a contradiction, represent a misinterpretation of Kani's example.\par
The following conjecture (which is sometimes referred to as Iitaka-Severi Conjecture) represents a generalization of Theorem \ref{deFSev}. As we shall see in the proof of Theorem \ref{section3thm}, the
difficulty in proving it lies in the fact (and only in the fact) that there
are no uniform birational embedding theorems for manifolds with big
canonical line bundle in higher dimensions (not even ineffective ones).
\begin{conjecture}[Iitaka-Severi]
Let $X$ be a compact complex manifold. Then the number of birational
equivalence classes of compact complex manifolds of general type having a
member $Y$ for which there exists a dominant rational map from $X$ to $Y$ is
finite.
\end{conjecture}
Our result in this section is the following effective generalization of Theorem \ref{deFSev}.
\begin{theorem}\label{section3thm}
Let $X$ be an $n$-dimensional compact complex manifold whose canonical line bundle $K_X$ is ample. Then the number $\mathcal F(X)$ of surjective holomorphic maps $f:X\to Y$, where $Y$ is any $n$-dimensional compact complex manifold with ample canonical bundle, is no more than
\begin{equation*}
2^nk^nK_X^n\cdot{(N+1)\cdot 2^nk^nK_X^n\choose
N}^{(N+1)(2^nk^nK_X^n{2^nk^nK_X^n+n-1\choose n}+{2^nk^nK_X^n+n-1\choose
n-1})},
\end{equation*}
where
$$N=(k^nK_X^n+n)^2-1,$$
and $k=k(n)$ is the effective very ampleness bound from {\rm Theorem \ref{effpluri}}.
\end{theorem}
\begin{remark} Although the bound for $\mathcal F(X)$ looks somewhat
complicated, its behavior with respect to $n$ and $K_X^n$ is easy to determine using Sterling's formula.
Namely, there exist explicit (exponential) functions $\alpha(n), \beta(n)$
such that
$$\mathcal F (X)\leq (\alpha(n)K_X^n)^{\beta(n)(K_X^n)^{n+5}}.$$
\end{remark}
In the case of $Y$ being a compact complex surface of general type, effective bounds in the same spirit as our Theorem \ref{section3thm} were given by Tsai in \cite{TsaiJAG} (see also \cite{TsaiIMRN}, \cite{TsaiCrelle} and also \cite{Guerra}). Ineffective results related to Theorem \ref{section3thm} have also been established in \cite{MDLM}, \cite{Maehara}, \cite{BD}.\par
It is well known that a result of the type of Theorem \ref{section3thm}
can be proved by what is commonly referred to as a ``boundedness
and rigidity argument''. The basic idea is to show that the
objects in question can be associated to Chow points in a finite (or even effectively finite) number of Chow varieties and that in any irreducible component the Chow points can correspond to at most one of the objects in question. Then, clearly, the number of the objects in question is no more than the number of relevant
irreducible components of the Chow varieties. Keeping this strategy in mind, we now start the
\begin{proof}[Proof of {\rm Theorem \ref{section3thm}}]
Let $f:X\to Y$ be one of the maps under consideration and let $\Gamma_f$
denote its graph. Let $p_1,p_2$ denote the two canonical projections of
$X\times Y$. Let $\phi_f$ denote the isomorphism $X\to \Gamma_f, x \mapsto
(x,f(x))$. The line bundle $p_1^*(kK_X)\otimes p_2^*(kK_Y)$ is
very ample on $X\times Y$ and embeds $\Gamma_f\subset X\times Y$ into
\begin{eqnarray*}
X\times \P^{h^0(Y,kK_Y)-1}\hookrightarrow\P^{h^0(X,kK_X)\cdot h^0(Y,kK_Y)-1}.
\end{eqnarray*}
Due to Proposition \ref{proph0bound} and the fact that $K_X^n\geq K_Y^n$ (see below), we can assume $\Gamma_f$ to be embedded into
\begin{equation*}
X\times\P^{N_1}\hookrightarrow P^N
\end{equation*}
with 
\begin{equation*}
N_1 := k^nK_X^n+n-1
\end{equation*}
and
\begin{equation*}
N := (k^nK_X^n+n)^2-1.
\end{equation*}\par
The degree of $\Gamma_f$ (measured in $\P^N$) can be estimated as follows:
\begin{eqnarray*}
\deg(\Gamma_f)&=&\int_{\Gamma_f}(p_1^*c_1(kK_X)+p_2^*c_1(kK_Y))^n\\
&=&\int_X\phi_f^*(p_1^*c_1(kK_X)+p_2^*c_1(kK_Y))^n\\
&=&\int_X (c_1(kK_X)+f^*c_1(kK_Y))^n\\
&\leq&\int_X (2c_1(kK_X))^n\\
&=&2^nk^nK_X^n.
\end{eqnarray*}
Note that the inequality is due to the fact that $K_X=f^*K_Y+D$, where $D$
is an effective divisor and $K_X$ and $f^*K_Y$ are ample, whence
\begin{eqnarray*}
K_X^{n-j}\left(f^*K_Y\right)^j&=&K_X^{n-j-1}\left(f^*K_Y+D\right)\left(f^*K_Y\right)^j\\
&=&K_X^{n-j-1}\left(f^*K_Y\right)^{j+1}+K_X^{n-j-1} . D .
\left(f^*K_Y\right)^j\\
&\geq& K_X^{n-j-1}\left(f^*K_Y\right)^{j+1}
\end{eqnarray*}
for $j=0,\ldots,n-1$. In particular, this computation yields
\begin{equation*}
K_X^n\geq (f^*K_Y)^n\geq K_Y^n,
\end{equation*}
which we used previously.\par
We now come to the rigidity part of the proof. In \cite[Corollary 3.2]{TsaiCrelle} it is shown that if $\pi:Z\to \Delta$ is a holomorphic family of smooth projective varieties of general type over a disk with $Z_0\cong Y$, there is no surjective holomorphic map $F:X\times\Delta\to Z$ with $F(X\times\{t\})=\pi^{-1}(t)$ unless $Z\cong Y\times \Delta$ and $F(\cdot,t)$ is independent of $t$.\par
Now take an irreducible component $I$ of $\Chow_{n,d}(X\times\P^{N_1})$ that contains a point corresponding to one of our graphs $\Gamma_f \subset X\times Y\subset X\times\P^{N_1}$. According to our previous boundedness considerations, we have $d\leq 2^nk^nK_X^n$. To be able to apply the rigidity property stated above, we need the following parametrization statement. For the details of its proof, we refer to \cite[Section 3]{Maehara}, noting that our situation is essentially the same as the one treated by Maehara.\par
There is a Zariski-open subset $U\subset I$ such that all Chow points $[\Gamma] \in U$ correspond to surjective holomorphic maps $f_{[\Gamma]}:X\to Y_{[\Gamma]}$ with $Y_{[\Gamma]}\subset\P^{N_1}$ being an $n$-dimensional projective manifold of general type. Moreover, $U$ contains all Chow points $[\Gamma_f]\in I$ that come from graphs of maps $f:X\to Y$ of the type considered in the statement of the Theorem.\par
Based on this parametrization statement, the above-mentioned rigidity property implies that for $[\Gamma_1],[\Gamma_2] \in U$, we have $f_{\Gamma_1}=f_{\Gamma_2}$, i.e. the number $\mathcal F(X)$ is no more than the number of relevant irreducible components of $\Chow_{n,d}(X\times\P^{N_1})$ for $d=1,\ldots, 2^nk^nK_X^n$. Clearly, only those components of $\Chow_{n,d}(X\times\P^{N_1})$ are relevant whose general points represent irreducible cycles.
However, from \cite{Kollarbook} (and also \cite{Guerra}), the following proposition is known.
\begin{proposition}[\cite{Kollarbook}, \cite{Guerra}]
Let $W\subset \P^n$ be a projective variety defined by equations of degree no more than $\tilde \delta$. Let $\Chow'_{k,\delta}(W)$ denote the union of those irreducible components of $\Chow_{k,\delta}(W)$ whose general points represent irreducible cycles. Then the number of irreducible components of $\Chow'_{k,\delta}(W)$ is no more than
\begin{equation*}
{(n+1)\max\{\delta,\tilde \delta\}\choose n}^{(n+1)(\delta {\delta+k-1\choose k }+{\delta+k-1 \choose k-1})}.
\end{equation*} 
\end{proposition}
\begin{remark}
Bounds on the complexity (i.e. the number of irreducible components) of Chow varieties have previously been produced by a number of authors. For example, the problem was extensively studied in Catanese's \cite{Ca}, and also in the papers of Green-Morrison (\cite{GM}) and Tsai (\cite{TsaiIMRN}). A new approach to handling Chow varieties of $1$-dimensional cycles is introduced in \cite{heiereffshaf}.
\end{remark}
Since the degrees of the defining equations of $X\times \P^{N_1}\subset \P^N$ under the Segre embedding are no more than $k^nK^n_X$, we conclude that our cardinality $\mathcal F(X)$ can be estimated from above by
\begin{eqnarray*}
&&\sum_{d=1}^{2^nk^nK_X^n} \# \text{ of irreducible components of }\Chow'_{n,d}(X\times\P^{N_1})\\
&\leq&2^nk^nK_X^n\cdot{({N}+1)\cdot 2^nk^nK_X^n\choose
{N}}^{({N}+1)(2^nk^nK_X^n{2^nk^nK_X^n+n-1\choose
n}+{2^nk^nK_X^n+n-1\choose n-1})}.
\end{eqnarray*}
\end{proof}\par
We remark that the nonequidimensional case (i.e. $\dim X > \dim Y$) can be reduced to our Theorem \ref{section3thm} by taking hyperplane sections. We shall express this fact as follows.
\begin{corollary}\label{nonequidimcase}
If we take the targets $Y$ in {\rm Theorem \ref{section3thm}} to be $n'$-dimensional with $n-n'>0$, then the analogous cardinality $\mathcal F_{n'}(X)$ is no more than the number obtained when replacing $K_X^n$ with $((n-n')k+1)^{n'}k^{(n-n')}K_X^n$ in the bound obtained in {\rm Theorem \ref{section3thm}}. 
\end{corollary}
\begin{proof}
We keep the notation from the proof of Theorem \ref{section3thm}. For a generic hyperplane section $X\cap H$ of $X$ in $\P^{N_1}$, the restriction of the maps in question to $X\cap H$ is still surjective (for an easy proof of this fact see \cite{MDLM}). Therefore, after taking $n-n'$ general hyperplane sections, we obtain an $n'$-dimensional submanifold $\tilde X$ to which we can apply Theorem \ref{section3thm}.\par
An $(n-n')$-fold iteration of the adjunction formula yields
\begin{eqnarray*}
K_{\tilde X}^{n'}&=&(\frac 1 k \mathcal O(1)|_{\tilde X}+(n-n')\mathcal O(1)|_{\tilde X})^{n'}\\
&=&((n-n')+\frac 1 k)^{n'}k^nK_X^n\\
&=&((n-n')k+1)^{n'}k^{(n-n')}K_X^n.
\end{eqnarray*} 
\end{proof}
\begin{remark}
The strategy of an effective boundedness and rigidity proof can be used in a number of similar settings. For example, in the paper \cite{heiereffshaf}, a uniform effective bound is established for the finiteness statement of the Shafarevich Conjecture over function fields (Theorem of Parshin-Arakelov). The arguments in that paper are more delicate due to the more complicated situation (one has to deal with moduli maps instead of maps of the form $f:X\to Y$), but the underlying principle is essentially the same.
\end{remark}
\begin{acknowledgement} It is a great pleasure to thank Professor Yum-Tong Siu for many invaluable discussions on (effective) algebraic geometry in general and finiteness theorems of the type discussed in the present paper in particular. These discussions took place while I enjoyed the generous hospitality of the Mathematics Department of Harvard University and the Institute of Mathematical Research at the University of Hong Kong. It is with sincere gratitude that I acknowledge support through the Schwerpunktprogramm ``Globale Methoden in der komplexen Geometrie'' of the Deutsche Forschungsgemeinschaft through the chair of Professor Alan Huckleberry at Bochum University. 
\end{acknowledgement}

\end{document}